\documentclass{amsart}

\usepackage{amssymb}

\theoremstyle{plain}
\newtheorem{theorem}{Theorem}[section]

\newtheorem{proposition}[theorem]{Proposition}

\numberwithin{equation}{section}

\theoremstyle{definition}

\newtheorem{example}[theorem]{Example}

\theoremstyle{remark}
\newtheorem{remark}[theorem]{Remark}

\newcommand{\R}{\mathbb{R}}
\newcommand{\C}{\mathbb{C}}
\newcommand{\Z}{\mathbb{Z}}

\newcommand{\Su}{\operatorname{SU}}
\newcommand{\So}{\operatorname{SO}}
\newcommand{\Gl}{\operatorname{GL}}
\newcommand{\su}{\mathfrak{su}}
\newcommand{\fg}{\mathfrak{g}}

\DeclareMathOperator {\CP} {{\mathbb C}{P}}
\DeclareMathOperator {\RP} {{\mathbb R}{P}}

\newcommand{\ow}{\omega}
\newcommand{\orbit}{\mathcal{O}}

\newcommand{\la}{\langle}
\newcommand{\ra}{\rangle}

\newcommand{\ma}[1]{\begin{pmatrix} #1 \end{pmatrix}}

\newcommand     {\comment}[1]   {}
\newcommand{\mute}[2] {}
\newcommand     {\printname}[1] {}

\newcommand{\labell}[1] {\label{#1}\printname{#1}}

\begin{document}

\title{Nonstandard Lagrangian submanifolds in $\CP^n$}
\author{River Chiang}
\begin{abstract}
We use Hamiltonian actions to construct nonstandard 
(as opposed to $\RP^n$ and $T^n$)
Lagrangian submanifolds of $\CP^n$. First of all, a 
quotient of $\RP^3$ by the dihedral group $D_3$ is a 
Lagrangian submanifold of $\CP^3$. Secondly, $\Su(n)/\Z_n$ 
are Lagrangian submanifolds of $\CP^{n^2-1}$.
\end{abstract}
\maketitle

\comment{only need isotropy orbit with the right dimension}

\comment{$S^3=\So(4)/\So(3)$ in $M^6$?}

\section{Introduction}

Let $(M, \ow)$ be a symplectic manifold. A submanifold $L \subset M$
is \textbf{Lagrangian} if $\dim L = \frac{1}{2}\dim M$ and 
$\ow$ vanishes on $T(L)$.

Using Floer homology, Paul Biran has obtained information on the
topology of Lagrangian submanifolds of some closed manifolds. 
One of Paul Biran's theorems concerns the Lagrangian submanifolds of
$\CP^n$ \cite{B-1, B-2}: 

\begin{theorem}\labell{T:b1}
Let $L$ be a Lagrangian submanifold of $\CP^n$ such that 
$H_1(L;\Z)$ is $2$-torsion (that is, $2H_1(L;\Z)=0$). Then
\begin{enumerate}
\item $H^*(L; \Z_2) \cong H^*(\RP^n; \Z_2)$ as graded vector spaces. 
\item If $n$ is even, the isomorphism in $(1)$ is an isomorphism 
of graded algebras. 
\end{enumerate} 
\end{theorem}

Seidel \cite{Se} also has a theorem with regard to the same subject:

\begin{theorem}\labell{T:se}
Let $L$ be a Lagrangian submanifold of $\CP^n$. Then
\begin{enumerate}
\item $H^1(L;\Z_{2n+2})\neq 0$.
\item $H^1(L;\Z_{2n+2})\ncong (\Z_2)^g$ for any $g \geq 2$.
\item If $H^1(L;\Z_{2n+2}) \cong \Z_2$, then $H^i(L;\Z_2) \cong \Z_2$ 
for all $i=0, \dots, n$.
\end{enumerate}
\end{theorem}

To test the limits of these topological restriction on Lagrangian submanifolds 
of $\CP^n$, we construct new examples using Hamiltonian actions. The 
constructions in this paper relies on the following well-known fact:

\mute{}{By starting with the standard Lagrangian embeddings of $\RP^n$ and
$T^n$ and applying Lagrangian surgery, one can construct many
different Lagrangian submanifolds of $\CP^n$.} 

\comment{does surgery give this example? does surgery give examples such that 
$H^*(L;\Z_2)$ is not $H^*(\RP^n; \Z_2)$? It seems that $H_1$ is the only 
important thing}

\begin{proposition}\labell{P:iso}
Let a compact Lie group $G$ act on a symplectic manifold $(M, \ow)$ 
in a Hamiltonian fashion with a moment map 
$\Phi\colon M \to \fg^*$. Assume $\alpha \in \fg^*$ is a fixed point of the
coadjoint action. Then for any $m\in \Phi^{-1}(\alpha)$,
the orbit $G\cdot m$ is an isotropic submanifold of $M$. 
\end{proposition}

\begin{proof}
An orbit is smooth. We only need to show that it is isotropic. 
Let $\xi$ and $\eta$ be any vector fields in the Lie algebra $\fg$ of $G$ and
let $\xi_M$ and $\eta_M$ be their induced vector fields on $M$.  
\[ 
\begin{split}
\ow(m)(\xi_M(m),\eta_M(m)) 
&= d\la \Phi(m), \xi \ra (\eta_M(m)) \\
&=\la d\Phi(\eta_M), \xi \ra \\
&=\la \eta_{\fg^*}(\Phi(m)), \xi \ra\\
&=\la \frac{d}{dt}\bigr|_0\exp t\eta \cdot \alpha, \xi \ra\\
&=\la 0, \xi \ra 
\end{split}
\]
Therefore $\ow(m) = 0$ for all $m\in \Phi^{-1}(\alpha)$. In particular, 
the orbit $G\cdot m$ is a smooth isotropic submanifold of $M$.  
\end{proof}

In particular, any orbit of an abelian action is isotropic as well as
any orbit in the zero level set of a nonabelian action.

\begin{example}
Consider the action of the torus $T^n$ on $\CP^n$ given by
\[
(\lambda_1, \dots, \lambda_n)\cdot[z_0, z_1, \dots, z_n]
= [z_0, \lambda_1^{-1}z_1, \dots,\lambda_n^{-1}z_n].
\] This action is Hamiltonian and its moment map is 
\[
\Phi([z_0, z_1, \dots, z_n])=\frac{1}{2}\left(\frac{|z_1|^2}{\|z\|^2}, \dots, 
\frac{|z_n|^2}{\|z\|^2}\right),
\]where $\|z\|^2 = \sum_{i=0}^n|z_i|^2$. The free orbit $T^n$ is a 
Lagrangian submanifold of $\CP^n$.  
\end{example}

Our first example is a Lagrangian submanifold of $\CP^3$. 
Its topology 
has an interesting deviation from Theorems \ref{T:b1}
and \ref{T:se}. 

\begin{theorem}
There exists a Lagrangian submanifold of $\CP^3$, which is a quotient of 
$\RP^3$ by the dihedral group $D_3$.
\end{theorem}

We discuss the topology of this submanifold in 
Section \ref{S:top}.     

We use different actions to construct a second family of examples.
It provides a generalization of the standard 
Lagrangian submanifold $\RP^3$ of $\CP^3$. 

\begin{theorem}
$\Su(n)/\Z_n$ is a Lagrangian submanifold of $\CP^{n^2-1}$.
\end{theorem}

\section{First example}\labell{S:con1}

View $\C^4$ as the space of homogeneous polynomials 
of degree $3$ in two variables. We can define the
irreducible $\Su(2)$ representation $\rho\colon \Su(2) \to \Gl(\C^4)$ by
\[
\rho(A)(p)(z) = p(zA)
\]
where $p \in \C^4$, $A \in \Su(2)$ and $z = (x, y) \in \C^2$.  
The moment map for this representation as a map from $\C^4$ into $\R
\times \C \cong \su(2)^*$ is given by
\[
(u_0, u_1, u_2, u_3) \mapsto
(\frac{3}{2}|u_0|^2+\frac{1}{2}|u_1|^2-\frac{1}{2}|u_2|^2-\frac{3}{2}|u_3|^2,
u_0\overline{u_1}+u_1\overline{u_2}+u_2\overline{u_3}). 
\]
So the zero level set of the moment map must satisfy
\[
\begin{cases}
&3|u_0|^2+|u_1|^2=|u_2|^2+3|u_3|^2,\\
&u_0\overline{u_1}+u_1\overline{u_2}+u_2\overline{u_3}=0.
\end{cases}
\]

This $\Su(2)$ action commutes with the diagonal $S^1$ action on
$\C^4$, and therefore descends to an action on $\CP^3$. 

Clearly the point $[1, 0, 0, 1]$ is in the zero level set. 
We know that  
\[
(x, y)\begin{pmatrix}
\alpha & \beta\\
-\overline{\beta}&\overline{\alpha}
\end{pmatrix} = (\alpha x- \overline{\beta}y, \beta x + \overline{\alpha}y),
\] where $|\alpha|^2+|\beta|^2=1$. To compute the stabilizer, we set 
$\lambda (x^3+y^3) = (\alpha x- \overline{\beta}y)^3 + (\beta x +
\overline{\alpha}y)^3$ for any possible $\lambda$ in $S^1$. 
Direct computation shows that $\lambda = \pm 1$ when $\beta =0$ and
$\lambda = \pm i$ when $\alpha = 0$. So the stabilizer is 
generated by 
\begin{equation}\label{E:generator}
\begin{pmatrix}
0 & i \\
i & 0
\end{pmatrix} \quad \text{and} \quad
\left\{
\begin{pmatrix}
\alpha & 0 \\
0 &\overline{\alpha}
\end{pmatrix}, \alpha^6=1 \right\}
\end{equation}

Note this $\Su(2)$ action is not effective on $\CP^3$. It induces
an effective $\So(3)$ action and the stabilizer for $[1,0,
0,1]$ is the semidirect product $\Z_2 \ltimes \Z_3$, namely, the 
dihedral group $D_3$. 

Therefore, we have a Lagrangian $L = \So(3)/{D_3} \cong \RP^3/{D_3}$ 
in $\CP^3$. 

In fact, this is the only orbit in the zero level set. 
According to the Guillemin-Marle-Sternberg local 
normal form theorem, a neighborhood of the orbit
of $[1,0,0,1]$ is equivariantly symplectomorphic to a neighborhood of
the zero section in $Y = \So(3)\times_{D_3} \R^3$. It is easy to check that
there exists one single orbit in the zero level set in $Y$. Since the zero 
level set is connected and since the orbit is closed, it follows that the
zero level set in $\CP^3$ consists of one orbit only.

\section{Topology for the first example}\labell{S:top}

Let $\Gamma$ denote the group generated by the elements in 
(\ref{E:generator}). Then the fundamental group of $L$ is
$\pi_1(L) = \Gamma$ since $L= \RP^3/D_3 = S^3/\Gamma$.

Let $b$ denote $\left(\begin{smallmatrix}
0 & i \\
i & 0
\end{smallmatrix}\right)$ and $a$ denote 
$\left(
\begin{smallmatrix}
\alpha & 0 \\
0 &\overline{\alpha}
\end{smallmatrix}\right)$, $\alpha^6=1$. Then $b^2 = -1$ and $ba = a^5b$. 
The commutator subgroup $A = \{\,x^{-1}y^{-1}xy\,|\,x, y \in 
\pi_1(L)\, \}$ is $\{\,1, a^2, a^4\,\}$. By Hurewicz Theorem, 
the first homology group $H_1(L; \Z)$ is $\{\,1, b, b^2, b^3\,\}
\cong \Z_4$.  

Using the Universal Coefficient Theorem and Poincar\'e duality, we
can calculate that
\[
\begin{aligned}
H_0(L;\Z)&=\Z,\qquad &H^0(L;\Z)&=\Z, \qquad &H^0(L;\Z_2)=\Z_2,\\
H_1(L;\Z)&=\Z_4,\qquad &H^1(L;\Z)&=0, \qquad &H^1(L;\Z_2)=\Z_2,\\
H_2(L;\Z)&=0, \qquad &H^2(L;\Z)&=\Z_4, \qquad &H^2(L;\Z_2)=\Z_2,\\
H_3(L;\Z)&=\Z, \qquad &H^3(L;\Z)&=\Z, \qquad &H^3(L;\Z_2)=\Z_2.
\end{aligned}
\]
Also, $H^1(L;\Z_8) = \Z_4$.

The Lagrangian $L$ we found does not 
satisfy the assumption of $2$-torsion as in Theorem \ref{T:b1} and 
Theorem \ref{T:se}(3), 
but $H^*(L; \Z_2)$ is isomorphic to 
$H^*(\RP^3; \Z_2)$ as graded vector spaces. 

\comment{can we say something about generators?}

\comment{add a general statement about extending the theorems?}

\section{Second family of examples$^1$}\labell{S:con2}

\addtocounter{footnote}{+1}
\footnotetext{I would like to thank Mich\`ele Audin for showing me a
beautiful alternative construction of these examples in terms of 
real forms of Lie groups \cite{Au}.}

Consider $\C^{n^2}$ as the set of $n$ by $n$ matrices. $\Su(n)$ acts 
naturally by left multiplication. This action is Hamiltonian and its 
moment map is given by
\[
\Phi(Z) = \frac{i}{2}\,ZZ^* + \frac{1}{2i}\,\text{I}
\] where $I$ is the identity matrix. 

This $\Su(n)$ action 
commutes with the circle action of multiplication by $\lambda\,\text{I}$.
Therefore, it descends to an action on $\CP^{n^2-1}$ where the center $\Z_n$
of $\Su(n)$ acts trivially. Therefore, 
a free orbit in the zero level set of $\C^{n^2}$ descends to an orbit in 
the zero level set of $\CP^{n^2-1}$ with stabilizer $\Z_n$. Since
the dimension of $\Su(n)$ is $n^2-1$, by Proposition
\ref{P:iso}, $\Su(n)/\Z_n$ is a Lagrangian submanifold of $\CP^{n^2-1}$. 

For $n=2$, the above-mentioned orbit $\Su(2)/\Z_2$ is the 
standard Lagrangian $\RP^3$ in $\CP^3$. Indeed it is the projective
space of the real $4$ dimensional space $V$ spanned by the following
matrices: 
\[
\ma{1&0\\0&-1}, \ma{i&0\\0&i}, \ma{0&1\\1&0}, \ma{0&i\\-i&0}.
\]


\begin{thebibliography}{9999999}
   
  \bibitem[Au]{Au} M. Audin, Unpublished notes. 

  \bibitem[A-L-P]{A-L-P} M. Audin, F. Lalonde and L. Polterovich, {\em
     Symplectic rigidity: Lagrangian submanifolds.} In Holomorphic
   curves in symplectic geometry.  Edited by M. Audin and J.
   Lafontaine. Progress in Mathematics, \textbf{117}. Birkh\"{a}user
   Verlag, Basel, 1994.

  \bibitem[B-1]{B-1} P. Biran, {\em Homological uniqueness of 
     Lagrangian submanifolds.} In preparation.

  \bibitem[B-2]{B-2} P. Biran, {\em Geometry of symplectic intersections.} 
	 Proceedings of the International Congress of Mathematicians 
	(Beijing 2002), Vol. II, 241--255.

   
   
   
   
   
   
     
   
   
   
  \bibitem[Se]{Se} P. Seidel, {\em Graded Lagrangian submanifolds.}
   Bull. Soc. Math. France \textbf{128} (2000), no. 1, 103--149.

\end{thebibliography}
\end{document}